\documentclass[a4paper,11pt]{article}
\usepackage{amsfonts}
\usepackage{pstricks}
\usepackage{pst-node,pst-tree,pst-eps,xspace}
\IfFileExists{pst-beta.tex}
             {\input{pst-beta.tex}}
             {\RequirePackage{pst-tree}}

\usepackage{latexsym}
\usepackage{amsthm}
\usepackage{fancyheadings}
\usepackage{amsmath}
\usepackage{graphicx}
\usepackage{float}
\usepackage{multirow}
\usepackage{array}
\usepackage{version}
\usepackage[font=small]{caption}
\restylefloat{table}
\usepackage{mathtools}

\newtheorem{Co}{Corollary}
\newtheorem{Lem}{Lemma}

\newtheorem{The}{Theorem}
\newtheorem{Pro}{Proposition}
\theoremstyle{definition}
\newtheorem{De}{Definition}
\newtheorem{Rem}{Remark}
\newtheorem{Exam}{Example}

\newcommand{\bsb}[1]{\boldsymbol{#1}}

\newcommand{\ccoomm}[1]{#1}

\author{
A. Bernini\thanks{Dipartimento di Matematica e Informatica ``Ulisse Dini", Universit\`a
degli Studi di Firenze, Viale G.B. Morgagni 65, 50134 Firenze, Italy.
{\tt \ \{bernini\}\{bilotta\}\{pinzani\}@dsi.unifi.it
}} 
\and S. Bilotta$^*$\and R. Pinzani$^*$
\and A. Sabri\thanks{LE2I, Universit\'e de Bourgogne, BP 47 870, 21078 Dijon Cedex, France.  
{\tt \{ahmad.sabri\}\{vvajnov\}@u-bourgogne.fr }}
\and V. Vajnovszki${^\dag}$
}
\title{Gray code orders for $q$-ary words avoiding a given factor}
\begin{document}
\maketitle

\begin{abstract}
Based on BRGC inspired order relations we define Gray codes and give a generating 
algorithm for $q$-ary words avoiding a prescribed factor. 
These generalize an early 2001 result and a very recent 
one published by some of the present authors, and can be seen as 
an alternative to those of Squire
published in 1996. Among the involved tools, we make use of
generalized BRGC order relations, ultimate periodicity
of infinite words, and word matching techniques.
 
\end{abstract}

\section{Introduction}

A very special way for listing a class of combinatorial objects
is the so called \emph{combinatorial Gray code}, where two consecutive objects differ 
`in some pre-specified small way'~\cite{JWW}.
In \cite{Wa1} a general definition is given, where 
a Gray code is defined as
`an infinite set of word-lists
with unbounded word-length such that the Hamming distance between any two
adjacent words is bounded independently of the word-length' (the Hamming
distance is the number of positions in which the words differ).

In \cite{Guibas_odlyz} Guibas and Odlyzko enumerated the set of length $n$ words 
avoiding an arbitrary factor, and a systematic construction and enumeration results
for particular factor avoidance in binary case are considered in \cite{BPP_2012,BMPP_2012}.
In Gray code context, Squire in his early paper \cite{Squ_96} explores 
the possibility of listing factor avoiding words such that  
consecutive words differ in only one position, and by $1$ or $-1$ in this position,
and in \cite{Vaj_2001} is given a Gray code and a generating algorithm
for binary words avoiding $\ell$ consecutive~1s. 
The result in \cite{Vaj_2001} was recently generalized in \cite{BBPV_2013}
where two Gray codes (one prefix partitioned and the other trace partitioned)
for $q$-ary words avoiding a factor constituted by $\ell$ consecutive
equal symbols are given.


Here, we adopt a different approach by relaxing
Squire's `one position constraint' and give Gray codes for length $n$ words avoiding 
any given factor, where consecutive words differ in at most three positions.
Our definitions for these Gray codes are based on two order relations 
inspired from the original {\it Binary Reflected Gray Code} \cite{Gray};
similar techniques were used previously (less or more explicitly) for other 
combinatorial classes, see for example 
\cite{Walsh_2000,Sab_Vaj_2013,Vaj_Ver_2011}
and the references therein.
More precisely, we characterize forbidden factors inducing 
{\it zero periodicity} (defined later), which is a crucial notion for our construction; 
and we show that the zero periodicity property of a forbidden factor 
is a sufficient condition for the set of words avoiding this factor 
when listed in the appropriate order to be a Gray code. 
However, this is not a necessary condition and we show that 
there are forbidden factors with no zero periodicity
property, and the set of words avoiding one of them when listed in the appropriate order 
yields a Gray code.
Also, among all $q^\ell$ forbidden factors of length $\ell$ on a $q$-ary alphabet,
all but $\ell+1$ of them induce zero periodicity; and
when a Gray code is prohibited by lack of zero periodicity property of the forbidden 
factor, we give a simple transformation of this factor which allows to eventually 
obtain the desired Gray code. 
Finally, we give a constant average time generating algorithm for these
Gray codes using $\ell \cdot q$ extra space and a 
Knuth-Morris-Pratt word matching technique \cite{Knu_Mor_Pra_1977}.
A {\tt C} implementation of the obtained  
algorithm is on the web site of the last author \cite{Vaj_web}.

Although in \cite{Li_Sawada} it is proved that the set of 
words avoiding a given factor is `reflectable'
under some conditions on the alphabet cardinality and the forbidden factor,
our construction yields Gray codes 
for any alphabet and forbidden factor, and has a natural algorithmic implementation.

\section{Notations and definitions}
\label{sec:notation}

\subsection*{Words over a finite alphabet}
An {\it alphabet} $A$ is simply a set of symbols,
and a length  $n$ {\it word} is a 
function $\{1,2,\ldots, n\}\rightarrow A$,  and $\epsilon$ is the 
empty (i.e., length zero) word.
We adopt the convention that lower case bold letters 
represent words, for example $\bsb a=a_1a_2\ldots a_n$; and
$A^n$ denotes the set of words of length $n$ over $A$,  
$A^*=\cup_{n\geq 0}A^n$ and 
$A^+=\cup_{n\geq 1}A^n$. 
For ${\bsb a}\in A^*$, $|{\bsb a}|$ denotes the length of $\bsb a$, or equivalently,
the number of symbols in ${\bsb a}$,  and $|{\bsb a}|_{\neq 0}$ the number of non-zero 
symbols in ${\bsb a}$.
%
An infinite word is a function $\mathbb{N}\rightarrow A$, and $A^\infty$ is the set of infinite 
words over $A$.  
For ${\bsb a}\in A^*$,  $i\geq 0$, $\bsb a^i$ is the word obtained by $i$ repetitions of $\bsb a$
($\bsb a^0$ being the empty word $\epsilon$) and $\bsb a^\infty$ is the infinite
{\it periodic} word ${\bsb a}{\bsb a}{\bsb a}\ldots$.
The word ${\bsb a}\in A^\infty$ is {\it ultimately periodic} if there are 
${\bsb b}\in A^*$ and ${\bsb c}\in A^+$ 
such that ${\bsb a}={\bsb b}{\bsb c}^\infty$, and we say that $\bsb a$ has {\it ultimate period} $\bsb c$.
%
Incidentally, we will make use of {\it left infinite words}, which are
infinite words $\bsb a$ of the form $\bsb a=\ldots a_{-3}a_{-2}a_{-1}$.
Left infinite words are reverse of infinite words, and
formally a left infinite word is a function 
$\{\ldots, -3,-2,-1\}\rightarrow A$, and
for $\bsb a\in A^*$, $\bsb a^{-\infty}$ is the left infinite word 
$\ldots \bsb a\bsb a\bsb a$.

The word $\bsb f$ is a {\it factor} of the word $\bsb a$ if there are words
$\bsb b$ and $\bsb c$ such that ${\bsb a}={\bsb b}{\bsb f}{\bsb c}$; 
when ${\bsb b}=\epsilon$ (resp. ${\bsb c}=\epsilon$), then ${\bsb f}$ is a 
{\it prefix} (resp.  {\it suffix})
of  ${\bsb a}$; and in this case, the prefix or the suffix is {\it proper} if 
$\bsb f\neq \bsb a$ and 
$\bsb f\neq \epsilon$.

For a word $\bsb a$ and a set of words $X$ we denote by $\bsb a|X$ the set of 
words in $X$ having the prefix $\bsb a$, and by $X(\bsb a)$ those avoiding $\bsb a$, i.e, 
the words in $X$ which do not contain $\bsb a$ as a factor.
Thus, for example, $\bsb p|X(\bsb f)$ is the set of words in $X$
having prefix $\bsb p$ and avoiding $\bsb f$.
Clearly $A^*({\bsb a})=\cup_{n\geq 0}A^n({\bsb a})$.

Through this paper we consider the alphabet $A_q=\{0,1,\ldots,q-1\}$
with $q\geq 2$.

\subsection*{Gray codes}
We will adopt the following definition:
a list of same length words is a {Gray code} if there is a $d$
such that the Hamming distance between any consecutive words in the list 
is bounded from above by $d$;
and often we refer to this list as a $d$-Gray code, and so, for example, a  $3$-Gray code is also a 
$4$-Gray  code.
In addition, if for any two consecutive words in the list the 
leftmost and the rightmost positions where they differ are 
separated by at most $e-1$ symbols, then the Gray code is 
called {\it $e$-close}.

\subsection*{Order relations}
Our constructions of Gray codes for factor avoiding words are based on two 
order relations on $A_q^n$ we will define below.
The first one captures the order induced by {\it $q$-ary Reflected Gray Code} 
\cite{Er}, which is the natural extension of {\it Binary Reflected Gray Code} 
introduced by Frank Gray \cite{Gray}; and the second one is a variation of 
the previous one.
\begin{De}[\cite{Vaj_2001,Vaj_Ver_2011,Sab_Vaj_2013}]
\label{De_RGCo_e}
Let $\bsb{s}=s_1s_2\ldots s_n$ and $\bsb{t}=t_1t_2\ldots t_n$
be two words in $A^n_q$, $k$ be the leftmost position where 
they differ, and $u=\sum_{i=1}^{k-1} s_i=\sum_{i=1}^{k-1} t_i.$
We say that $\bsb{s}$ is less than $\bsb{t}$ in 
{\it Reflected Gray Code order}, denoted by $\bsb{s}\prec\bsb{t}$,
if either
\begin{itemize}
\item $u$ is even and $s_k<t_k$, or
\item $u$ is odd and $s_k>t_k$.
\end{itemize}
\end{De}

It follows that the set $A_q^n$ listed in $\prec$ order
yields precisely the $q$-ary Reflected Gray Code (see \cite{Er,Will_85}) where two consecutive words 
differ in one position and by $1$ or $-1$ in this position.

For a set of same length words $X$, we refer to {\it $\prec$-first} 
(resp.{\it  $\prec$-last}) word in $X$ for the first (resp. last) word 
in $X$ with respect to $\prec$ order.

\medskip
\noindent
In the binary case, Definition \ref{De_RGCo_e} can be re-expressed as
following.
For $\bsb{s}$, $\bsb{t}$ and $k$ as in Definition \ref{De_RGCo_e},
let $v$ be the number of non-zero symbols in the 
length $k-1$ prefix of $\bsb{s}$ and of $\bsb{t}$.
Then $\bsb{s}\prec\bsb{t}$ if either
\begin{itemize}
\item
  $v$ is even and $s_k<t_k$, or
\item
  $v$ is odd and $s_k>t_k$.
\end{itemize}

By `adding' $u$ in Definition \ref{De_RGCo_e} and $v$ defined above 
we obtain a new order relation.

\begin{De}
\label{De_RGCo_O}
Let $\bsb{s}=s_1s_2\ldots s_n$ and $\bsb{t}=t_1t_2\ldots t_n$
be two words in $A^n_q$, $k$ be the leftmost position where 
they differ, $u=\sum_{i=1}^{k-1} s_i=\sum_{i=1}^{k-1} t_i$,
and $v$ be the number of non-zero symbols in the 
length $k-1$ prefix of $\bsb{s}$. 
We say that $\bsb{s}$ is less than $\bsb{t}$ in 
{\it Dual Reflected Gray Code order}, denoted by $\bsb{s}\triangleleft\bsb{t}$,
if either
\begin{itemize}
\item $u+v$ is even and $s_k<t_k$, or
\item $u+v$ is odd and $s_k>t_k$.
\end{itemize}
\end{De}

Clearly, listing a set of words in $\prec$ or in $\triangleleft$ order gives a prefix 
partitioned list, in the sense that words with the same prefix are consecutive.
See Table \ref{ex_A3_4} in Appendix for the set $A_3^4$ listed
in $\triangleleft$ order.

For a set of same length words $X$, we refer to {\it $\triangleleft$-first} 
(resp.{\it  $\triangleleft$-last}) word in $X$ for the first (resp. last) word 
in $X$ with respect to $\triangleleft$ order.
In the following, without explicitly precise otherwise 
we will consider $\prec$ order on $A_q^*$ when $q$ is even and $\triangleleft$
order when $q$ is odd.

The next remark says that $\triangleleft$ order yields a Gray code when $q$ is odd,
but generally for $q$ even, $A_q^n$ listed in $\triangleleft$ order is not a 
Gray code.

\begin{Rem}
\label{G_VoRGCO}
For any odd $q\geq 3$ and $n\geq 1$, the set $A_q^n$ listed in $\triangleleft$ 
order is a Gray code where two consecutive words differ in at most two
adjacent  positions. In addition, if ${\bsb s}$ and ${\bsb t}$ are 
two consecutive words in this list and $k$ is the leftmost position where
they differ, then 
\begin{itemize}
\item $s_k=t_k\pm 1$, and 
\item if $s_{k+1}\neq t_{k+1}$, then $\{s_{k+1},t_{k+1}\}=\{0,q-1\}.$
\end{itemize}

\end{Rem}

We define the {\it parity} of the word $\bsb a=a_1a_2\ldots a_n\in A_q^n$ according to two cases:
\begin{itemize}
\item when $q$ is even, the parity of $\bsb a$ is the parity of the integer
      $\sum_{i=1}^n a_i$, and
\item when $q$ is odd, the parity of $\bsb a$ is the parity of the integer
      $\sum_{i=1}^n a_i+|\bsb a|_{\neq0}$.
\end{itemize}

For example, if $\bsb a=0222$, then
\begin{itemize}
\item considering $\bsb a\in A_4^*$, the parity of $\bsb a$ is even, and given by 
      $0+2+2+2=6$, and
\item considering $\bsb a\in A_5^*$, the parity of $\bsb a$ is odd, and given by 
      $0+2+2+2+3=9$.
\end{itemize}

Now we introduce a critical concept for our purposes:
we say that the forbidden factor $\bsb f\in A_q^*$
{\it induces zero periodicity} on $A_q^\infty$ if for any $\bsb p\in A_q^*({\bsb f})$
the first and the last (as mentioned previously, with respect to 
$\prec$ order for even $q$, or $\triangleleft$ 
order for odd $q$) words in $\bsb p|A^\infty_q({\bsb f})$
both have ultimate period $0$. 
Consequently, if $\bsb f\in A_q^*$ {\it does not} induce zero periodicity on 
$A_q^\infty$, it follows that there exists a $\bsb p\in A_q^*(\bsb f)$ such that  
the first and/or the last word in $\bsb p|A_q^\infty(\bsb f)$
do not have ultimate period 0.

Whether or not $\bsb f$ induces zero periodicity on $A_q^\infty$
depends on $q$; and $q$ will often be understood from the context.
For example $\bsb f=3130$ induces zero periodicity on
$A_6^\infty$ but not on $A_4^\infty$.
Indeed, the first word in $\prec$ order in $A_6^\infty(\bsb f)$ and having prefix 
$313$ is $3135000\ldots$, and the last one is $3131500\ldots$;
whereas the last word in $\prec$ order 
in $A_4^\infty(\bsb f)$ and with the same prefix is the periodic word $31313\ldots$.
 
In Section \ref{Sec_Gray_codes} it is shown that the {\it Graycodeness} of the set $A_q^n(\bsb f)$
listed in the appropriate order is intimately related to that 
$\bsb f$ induces zero periodicity.

\subsection*{Outline of the paper}
Avoiding a factor of length one is equivalent to shrink the
underlying alphabet, but  for the sake of generality, 
we will consider most of the time forbidden factors of any positive length. 

In the next section we will characterize the forbidden factors $\bsb f$
inducing zero periodicity on $A_q^\infty$. 
By means of three sets $U_q,V,W_q\subset A_q^+$, these factors are 
characterized in Corollary \ref{Cor}.
The non-zero ultimate periods produced by forbidden factors that do not induce zero 
periodicity have the form $1(q-1)0^m$, $10^m$ or $(q-2)$,
see Propositions \ref{pro:end_0_even}, \ref{pro:end_0_odd} and \ref{pro:end_max}.

\medskip

Theorems \ref{the:general_end0_qEven} to \ref{the:general_max} and Proposition 
\ref{pro:general_other} in Section \ref{Sec_Gray_codes} prove that the property 
of $\bsb f$ to induce zero periodicity guarantees the set
$A^n_q({\bsb f})$ listed in the appropriate order
to be a Gray code. 
However,  this property of $\bsb f$ is not a necessary condition:
there are two `special' forbidden factors, namely $\bsb f=0^\ell$ and 
$\bsb f=(q-1)0^\ell$
belonging to $U_q$
(and so, which do not induce zero periodicity) but $A_q^n(\bsb f)$ listed in 
$\prec$ order is still a Gray code.  
These two cases are discussed in Section \ref{sec:particular}. Table \ref{summary}
summarizes the Graycodeness for the set $A_q^n(\bsb f)$ 
if $\bsb f$ does not induce zero periodicity. Section \ref{Sec_Gray_codes} 
ends by showing that simple transformations of forbidden factors $\bsb f$
which do not induce zero periodicity allow to obtain Gray code for the set 
$A_q^n(\bsb f)$.

\medskip

Finally, we present in Section \ref{sec:algorithm} an efficient generating algorithm 
for the obtained Gray codes.

\section{Periodicity}

As stated above, without restriction, the set $A_q^n$ 
listed in $\prec$ order is a $1$-Gray code for any $q\geq 2$, and
listed in $\triangleleft$ order a $2$-Gray code for odd $q\geq 3$.
Roughly, it is due to that for any $\bsb p\in A_q^*$ the
first and the last words---with respect to $\prec$ order for any $q\geq 2$,
or $\triangleleft$ order for odd $q\geq 3$---
in the set $\bsb p|A_q^\infty$ are among the three infinite words:
$\bsb p0^\infty$,  $\bsb p(q-1)0^\infty$,
and $\bsb p(q-1)^\infty$.
This phenomenon is no longer true if an arbitrary factor $\bsb f$ is forbidden.
For example, if $\bsb f=130$, then the $\prec$-last word in $03|A_4^\infty(\bsb f)$ is 
$0300000\ldots$, and the $\prec$-first one in $13|A_4^\infty(\bsb f)$ is $1313131\ldots$;
and  $0300000$ and $1313131$ are consecutive words in $A_4^7(\bsb f)$, in $\prec$ order.
Or, for $\bsb f=223$, the $\prec$-last word in $123|A_4^\infty(\bsb f)$
is $123300\ldots$ and the $\prec$-first one in $122|A_4^\infty(\bsb f)$ 
is $122222\ldots$; and $123300$ and $122222$ are consecutive words in $A_4^6(\bsb f)$, 
in $\prec$ order.

However, it is easy to understand the next remark.

\begin{Rem} 
\label{easy} 
If $\bsb f\in A_q^*$ is a forbidden factor ending by a symbol other
than $0$ or $q-1$, then for any $\bsb p\in A_q^*(\bsb f)$, 
$q\geq 2$ and even (resp. $q\geq 3$ and odd),
the set formed by the $\prec$-first and the $\prec$-last (resp. the $\triangleleft$-first and 
the $\triangleleft$-last) words in $\bsb p|A_q^\infty(\bsb f)$ is 
$\{\bsb p0^\infty,\bsb p(q-1)0^\infty\}$.
\end{Rem}

In other words, the previous remark says that any factor ending by a symbol other than
$0$ or $q-1$ induces zero periodicity. However, there exist forbidden 
factors ending by 0 or $q-1$ that do induce zero periodicity. For example,
with $\bsb f=120$, the $\prec$-first and
$\prec$-last words in 
$12|A_4^\infty(\bsb f)$ are $1230000\ldots$ and $1213000\ldots$.

In the following we will use (often implicitly) the next straightforward
remark which provides the form of the words on the right of a
fixed prefix $\bsb p\in A_q^*(\bsb f)$, with respect to the appropriate order.
It is obtained by the following observation: the first/last word in $\bsb p|A_q^n(\bsb f)$
is the appropriate prefix of the first/last word in $\bsb p|A_q^\infty(\bsb f)$.

\begin{Rem}
\label{easy_rem}
$ $
Let $q$ be even, $\bsb f\in A_q^*$ be a forbidden factor, 
$\bsb p,\bsb r \in A_q^*(\bsb f)$, and let $\bsb p$ have even (resp. odd) parity
such that  $\bsb{pr}\in A_q^*(\bsb f)$. 
Then:
\begin{itemize}
\item If $\bsb{pr}$ is a prefix of the $\prec$-first (resp. $\prec$-last) word in 
       $\bsb p| A_q^\infty(\bsb f)$, then $\bsb r$
       is the smallest word, in $\prec$ order, with this property; that is,
       if $\bsb s\in A_q^*$ with $|\bsb s|=|\bsb r|$ and $\bsb s\neq \bsb r$
       is such that $\bsb{ps}$
       is the prefix of some word in $\bsb p| A_q^\infty(\bsb f)$,
       then $\bsb r\prec \bsb s$.
\item If $\bsb{pr}$ is a prefix of the $\prec$-last (resp. $\prec$-first) word in 
       $\bsb p| A_q^\infty(\bsb f)$, then $\bsb r$
       is the largest word, in $\prec$ order, with this property; that is,
       if $\bsb s\in A_q^*$ with $|\bsb s|=|\bsb r|$ and $\bsb s\neq \bsb r$
       is such that $\bsb{ps}$
       is the prefix of some word in $\bsb p| A_q^\infty(\bsb f)$,
       then $\bsb s\prec \bsb r$.
\end{itemize}   

\noindent
Similar results hold when $q$ is odd by replacing $\prec$ by $\triangleleft$
and considering the words parity as in the definition given after 
Remark \ref{G_VoRGCO}.
\end{Rem}

Remark \ref{rem:category} below specifies the form of the 
first and last words in $A_q^n$, subject to the additional constraint that they do not 
begin by $0$ or $q-1$.
Later on we will see that when $\bsb f$ does not induce 
zero periodicity, then the possible non-zero periods of the first or 
the last word in $\bsb p|A_q^*(\bsb f)$ are related to those 
words. This remark will be used in the proofs of Propositions \ref{pro:end_0_even},
\ref{pro:end_0_odd} and~\ref{pro:end_max}.

\begin{Rem}
\label{rem:category}
$ $
\begin{itemize}
\item For $q$ even, the first word in $A_q^n$, with respect to  $\prec$ order, 
      which does not begin by $0$ is $1(q-1)0^{n-2}$.
\item For $q$ odd, the first word in $A_q^n$, with respect to  $\triangleleft$ 
      order, which does not begin by $0$ is $10^{n-1}$.
\item For $q\geq2$ (even or odd), the last word in $A_q^n$, with respect to the appropriate 
      order, which does not begin by a $q-1$ is $(q-2)(q-1)0^{n-2}$.
\end{itemize}
\end{Rem}

In the following we need the technical lemma below.

\begin{Lem}
\label{comb_w}
Let $\bsb u,\bsb g\in A_q^*$ and $\bsb v\in A_q^+$
be such that $\bsb g$ is a suffix of both $\bsb u$ and $\bsb u\bsb v$. 
Then there exist a $j\geq 0$ and a (possibly empty) suffix $\bsb w$ of $\bsb v$
such that $\bsb g=\bsb w\bsb v^j$
(or equivalently, $\bsb g$ is a suffix of the left infinite word  $\bsb v^{-\infty}$).
\end{Lem}
\ccoomm{
\proof
We prove the statement by induction on $k=\lfloor\frac{|\bsb g|}{|\bsb v|} \rfloor$.
When $k=0$, since the length of $\bsb g$ is less than that of $\bsb v$,
and $\bsb g$ is a suffix of $\bsb u\bsb v$ the statement follows by
considering $j=0$ and $\bsb w=\bsb g$.

Let now $k=\lfloor\frac{|\bsb g|}{|\bsb v|} \rfloor>0$.
In this case the length of $\bsb g$ is greater than that of $\bsb v$,
and it follows that $\bsb v$ is a suffix of both $\bsb g$ and $\bsb u$.
By considering $\bsb u'$ and $\bsb g'$ such that 
\begin{itemize}
\item $\bsb u=\bsb u'\bsb v$
\item $\bsb g=\bsb g'\bsb v$ 
\end{itemize}
we have that $\bsb g'$ is a suffix of both $\bsb u'$ and $\bsb u=\bsb u'\bsb v$.
Since $|\bsb g'|=|\bsb g|-|\bsb v|$ we have that 
$\lfloor\frac{|\bsb g'|}{|\bsb v|} \rfloor=k-1$ and the statement follows by induction
on $k$.
}
\endproof

\subsection{Forbidden factor ending by $0$ and not inducing zero periodicity}
%
We will determine, according to the parity of $q$,
the form of the first and the last words in 
$\bsb p|A_q^\infty(\bsb f)$ having no ultimate period $0$
for $\bsb f$ ending by $0$, and consequently the form of the forbidden factors
$\bsb f$ that do not induce 
zero periodicity.

\subsubsection*{$q$ even}
\medskip

\label{sec:q_even}
\begin{Pro}
\label{pro:end_0_even}
Let $q\ge2$ be even and $\bsb f\in A_q^+$ be a forbidden factor ending by $0$ 
and not inducing zero periodicity. Let also $\bsb p\in A_q^*(\bsb f)$ be such 
that one of the $\prec$-first or the $\prec$-last word in 
$\bsb p| A_q^\infty(\bsb f)$ does not have ultimate period $0$,
and let $\bsb a$ be this word.
Then $\bsb a$ is ultimately periodic, more precisely there
is an $m\geq0$ such that either
\begin{itemize}
\item[(i)]  $\bsb a=\bsb p 0^{i}(1(q-1)0^m)^\infty$, for some $i\leq m$, or
\item[(ii)] $\bsb a=\bsb p ((q-1)0^m1)^\infty$.
\end{itemize}
\end{Pro}
\ccoomm{
\proof
We will show that when $\bsb p$ has even (resp. odd) parity, then either
\begin{enumerate}
\item $\bsb a$ is the $\prec$-first (resp. $\prec$-last) word in 
      $\bsb p| A_q^\infty(\bsb f)$, and in this case $\bsb a$ has the form given 
      in point $(i)$ above, or
\item $\bsb a$ is the $\prec$-last (resp. $\prec$-first) word in 
      $\bsb p| A_q^\infty(\bsb f)$, and in this case $\bsb a$ has the form given
      in point $(ii)$ above.
\end{enumerate}
For the point 1, considering the parity of $\bsb p$ and
since  $\bsb a$ does not have ultimate period $0$, there is 
an $i\ge0$ such that $\bsb p 0^{i+1}$ contains the factor $\bsb f$, but $\bsb p 0^{i}$ 
does not. Now, by Remark \ref{rem:category},
since $\bsb f$ ends by a $0$, it follows that
there is an $m\ge 0$ such that $\bsb p 0^{i}\,1(q-1)0^m$ is a prefix of $\bsb a$,
but $\bsb p 0^{i}1(q-1)0^{m+1}$ is not. 
Thus, the length maximal $0$ suffix of $\bsb f$ is $0^{m+1}$, and reasoning in 
the same way, it follows that
there is an $m'\geq 0$ such that $\bsb p 0^{i}\,1(q-1)0^m\,1(q-1)0^{m'}$ is a prefix of 
$\bsb a$, but $\bsb p 0^{i}\,1(q-1)0^m\,1(q-1)0^{m'+1}$ is not.
Since $0^{m+1}$ is the length maximal $0$ suffix of $\bsb f$, necessarily $m'=m$, 
and the statement holds by iterating this construction.

Similarly, point 2 holds considering that $\bsb p(q-1)$ is a prefix of $\bsb a$ 
and there is an $m\geq 0$ 
such  that $\bsb p(q-1)0^{m+1}$ contains the factor $\bsb f$.
}
\endproof

Now we characterize the forbidden factors $\bsb f\in A_q^+$
ending by $0$, for even $q\geq2$, and not inducing zero periodicity.

For even $q\geq 2$, let define the set $U_q\subset A_q^+$ as
\begin{equation}
\label{eq:Uq}
U_q=\bigcup_{m\geq 0}\{
\bsb b0\, |\, \bsb b\ {\rm a\ suffix\ of\ }
(1(q-1)0^m)^{-\infty}\}.
\end{equation}

Alternatively, $U_q$ is the set of words of the form 
$\bsb b0$, where $\bsb b$ is either empty, or for some $m\geq 0$, a factor of $(1(q-1)0^m)^\infty$ ending by $0^m$ if $m>0$ and 
ending by $q-1$ elsewhere. Clearly,  $U_q$ contains exactly $n$
words of length $n$, for example, 
$U_4\cap A_4^5=\{00000,30000,13000,01300,13130\}$.

\begin{Pro}
\label{a_first_last1}
For even $q\ge 2$, if a forbidden factor $\bsb f\in A_q^+$ ending by $0$ does not 
induce zero periodicity, then $\bsb f\in U_q$.
\end{Pro}
\ccoomm{
\proof
If $\bsb f$ is such a factor, then by Proposition \ref{pro:end_0_even}
there is a $\bsb p\in A^*(\bsb f)$ and an $m\geq 0$ such that, $\bsb a$,
the $\prec$-first or the $\prec$-last word in $\bsb p| A_q^\infty(\bsb f)$
is 
\begin{itemize}
\item $\bsb a=\bsb p 0^{i}(1(q-1)0^m)^\infty$, for some $i\leq m$, or
\item $\bsb a=\bsb p ((q-1)0^m1)^\infty$.
\end{itemize}
Let $\bsb g$ be the word obtained from $ \bsb f$ after erasing its last $0$.
In the first case
it follows that $\bsb g$ is a suffix of both $\bsb p 0^{i}\,1(q-1)0^m$
and $\bsb p 0^{i}\,1(q-1)0^m\,1(q-1)0^m$, and by 
Lemma \ref{comb_w} the statement holds.
The proof is similar for the second case.
}
\endproof
%
%
%

\begin{Rem}
\label{f_b0_qeven}
If $\bsb f\in U_q$ and $q$ is even, then $\bsb f$ does not induce zero periodicity.
Indeed, let for example $\bsb f=\bsb b0$ with $\bsb b$ a suffix of $(1(q-1)0^m)^{-\infty}$ 
be as in relation (\ref{eq:Uq}). 
Then either the first word in $\bsb b| A_q^\infty (\bsb f)$ when 
$\bsb b$ has even parity, or the last word in $\bsb b| A_q^\infty (\bsb f)$ when 
$\bsb b$ has odd parity, has ultimate period $1(q-1)0^m$.
\end{Rem}

\begin{Exam}
Let $\bsb f=301300\in U_4$ be a forbidden factor and let consider the prefix 
$\bsb p=0021301\in A_4^*$. The $\prec$-first word in 
$\bsb p|A_4^\infty(\bsb f)$ is $\bsb p30(130)^\infty$.
\end{Exam}

Combining Proposition \ref{a_first_last1} and Remark \ref{f_b0_qeven}
we have the following theorem.

\begin{The}
\label{the:end_zero_even}
For even $q\geq2$, the forbidden factor $\bsb f\in A_q^+$ ending by $0$  
does not induce zero periodicity if and only if 
$\bsb f\in U_q$.
\end{The}

\subsubsection*{$q$ odd}
\medskip

Now we give the odd $q$ counterpart of the previous results.

\begin{Pro}
\label{pro:end_0_odd}
Let $q\geq3$ be odd and $\bsb f\in A_q^+$ be a forbidden factor ending by $0$ and 
not inducing zero periodicity. Let also $\bsb p\in A_q^*(\bsb f)$ be such that 
one of the $\triangleleft$-first or the $\triangleleft$-last word in $\bsb p| A_q^\infty(\bsb f)$
does not have ultimate period $0$, and let $\bsb a$ be this word.
Then $\bsb a$ is ultimately periodic, more precisely there
is an $m\geq0$ such that either
\begin{itemize}
\item $\bsb a=\bsb p 0^{i}(10^m)^\infty$, for some $i\leq m$, or
\item $\bsb a=\bsb p (q-1)(0^m1)^\infty$.
\end{itemize}
\end{Pro}
\ccoomm{
\proof
The proof is similar to that of Proposition \ref{pro:end_0_even} and considering the 
second point of Remark~\ref{rem:category}.
}
\endproof

Now we characterize the forbidden factor $\bsb f\in A_q^+$ ending by $0$, 
for odd $q\geq3$, and not inducing zero periodicity.

For $q\geq 3$, let define the set $V\subset A_q^+$ as
\begin{equation}
\label{eq:V}
V=\bigcup_{m\geq 0}\{
\bsb b0\, |\, \bsb b\ {\rm a\ suffix\ of\ }(10^m)^{-\infty}
\}.
\end{equation}

Alternatively, $V$ is the set of words of the form 
$\bsb b0$, where  $\bsb b$ is either empty, or   
for some $m\geq 0$, a factor of $(10^m)^\infty$ ending by $0^m$ if $m>0$ (and ending by $1$ elsewhere). 
Notice that $V$ does not depend on $q$, i.e. $V\subset A_q^+$ for any 
$q\geq 2$. 
Clearly,  $V$ contains exactly $n$ words of length $n$, for example,
$V\cap A_q^5=\{00000,10000,01000,10100,11110\}$, for any 
$q\geq2$.

Considering Proposition \ref{pro:end_0_odd} and the definition of $\triangleleft$
order relation, with the same arguments as in the 
proof of Proposition \ref{a_first_last1} we have the next result.
\begin{Pro}
\label{a_first_last}
For odd $q\ge 3$, if a forbidden factor $\bsb f\in A_q^+$ ending by $0$ does not 
induce zero periodicity, then $\bsb f\in V$.
\end{Pro}

\begin{Rem}
\label{f_b0_qodd}
If $\bsb f\in V$ and $q\geq 3$ is  odd, then $\bsb f$ does not induce zero periodicity
on  $A_q^\infty (\bsb f)$.
Indeed, let for example $\bsb f=\bsb b0$ with $\bsb b$ a suffix of $(10^m)^{-\infty}$
be as in relation (\ref{eq:V}). 
Then either the first word in $\bsb b| A_q^\infty (\bsb f)$ when 
$\bsb b$ has even parity, or the last word in $\bsb b| A_q^\infty (\bsb f)$ when 
$\bsb b$ has odd parity, has ultimate period $10^m$.
\end{Rem}

\begin{Exam}
Let $\bsb f=0100010000\in V$ a the forbidden factor and let consider the 
prefix $\bsb p=430100010\in A_5^*$. 
The $\triangleleft$-last  word in $\bsb p|A_5^\infty(\bsb f)$ 
is $\bsb p00(1000)^\infty$.
\end{Exam}

Combining Proposition \ref{a_first_last} and Remark \ref{f_b0_qodd}, we have the following 
theorem.
\begin{The}
\label{the:end_zero}
For odd $q\geq 3$, the forbidden factor $\bsb f\in A_q^+$ ending by $0$  
does not induce zero periodicity if and only if $\bsb f\in V$.
\end{The}

\subsection{Forbidden factor ending by $q-1$ and not inducing zero periodicity}
\medskip
The next proposition holds for $q\geq3$ (even or odd), and the case for $q=2$ is 
stated in the remark that follows it.

\begin{Pro}
\label{pro:end_max}
Let $q\ge3$ (even or odd) and $\bsb f\in A_q^+$ be a forbidden factor ending by 
$q-1$ and  not inducing zero periodicity. Let also $\bsb p\in A_q^*(\bsb f)$ be such that
one of the first or the last word in $\bsb p| A_q^\infty(\bsb f)$, with respect to the appropriate 
order, does not have ultimate period $0$, and let $\bsb a$ be this word.
Then $\bsb a=\bsb p(q-2)^\infty$.
\end{Pro}
\ccoomm{
\proof
Nor $\bsb p0$ neither $\bsb p(q-1)$  can not be a prefix of $\bsb a$; 
otherwise, in the first case $\bsb a=\bsb p0^\infty$ and in the second one 
$\bsb a=\bsb p(q-1)0^\infty$.
By the third point of Remark \ref{rem:category} and since $\bsb f$ ends by a $q-1$, it follows that 
$\bsb p(q-2)$ is a prefix of $\bsb a$,  but $\bsb p(q-2)(q-1)$ is not
(otherwise $\bsb a=\bsb p(q-2)(q-1)0^\infty$). Again, 
$\bsb p(q-2)(q-2)$ is a prefix of $\bsb a$, but $\bsb p(q-2)(q-2)(q-1)$ is not;
and finally $\bsb a=\bsb p(q-2)^\infty$.
}
\endproof

When $q=2$, the ultimate $(q-2)$ period of $\bsb a$ in Proposition \ref{pro:end_max}
becomes $0$ period, and so, for  $q=2$ any forbidden factor $\bsb f\in A_q^+$ ending by 
$q-1=1$ 
induces zero periodicity. Thus, below we will consider only factors ending by 
$q-1$ and not inducing zero periodicity only for $q\geq3$ (even or odd).

For $q\geq3$, let define the set $W_q$ as
\begin{equation}
\label{eq:Wq}
W_q=\bigcup_{\ell\geq 0}\{(q-2)^\ell(q-1)\}.
\end{equation}

With the previous terminology, $W_q$ is the set of words of the form
$\bsb b(q-1)$ with $\bsb b$ a suffix of $(q-2)^{-\infty}$.
Clearly, $W_q$ contains exactly one word of each length, and
for example, $W_4=\{3,23,223,2223,22223,\ldots\}$. 

\begin{Pro}
\label{a_first_last_2}
For $q\ge3$ (even or odd), if the forbidden factor $\bsb f\in A_q^+$ ending by 
$q-1$  does not induce zero periodicity,  then $\bsb f\in W_q$.
\end{Pro}
\ccoomm{
\proof
Let $\bsb f$ be such a factor, and $\bsb p\in A_q^*(\bsb f)$ such that, with respect to the appropriate order, 
the first  word in $\bsb p| A_q^\infty(\bsb f)$ has not ultimate period $0$ (the case of the first 
word being similar).
Let also $\bsb g$ be the (possibly empty) word obtained from $ \bsb f$ 
after erasing its last symbol $q-1$. By Proposition \ref{pro:end_max},
$\bsb g$ is a suffix of both $\bsb p(q-2)$ and $\bsb p(q-2)(q-2)$, and by 
Lemma \ref{comb_w} the statement holds.
}
\endproof

\begin{Rem}
\label{f_max}
If $\bsb f\in W_q$, then $\bsb f$ does not induce zero periodicity.
Indeed, let for example $\bsb b=(q-2)^\ell$, for some $\ell\geq 0$,
and $\bsb f=\bsb b(q-1)$ be as in relation (\ref{eq:Wq}).
Then the last word in 
$\bsb b| A_q^\infty (\bsb f)$ has ultimate period $(q-2)$.
\end{Rem}

\begin{Exam}
Let $\bsb f=223\in W_4$ be a forbidden factor and let consider the prefix 
$\bsb p=2322\in A_4^*$. The $\prec$-first word in $\bsb p|A_4^\infty(\bsb f)$ is $\bsb
p2^\infty$.
And when $\bsb f=12\in W_3$ and $\bsb p=01\in A_3^*$, 
the $\triangleleft$-last word in $\bsb p|A_3^\infty(\bsb f)=\bsb p1^\infty$.
\end{Exam}

Combining Proposition \ref{a_first_last_2} and Remark \ref{f_max}, we have the 
following theorem. 

\begin{The}
\label{the:end_max}
For $q\ge 3$ (even or odd), the forbidden factor $\bsb f\in A_q^+$ ending by 
$q-1$ does not induce zero periodicity if and only if $\bsb f\in W_q$.

\end{The}

Even we will not make use later, it is worth to mention the following remark.
\begin{Rem}
\label{useless}
For $q$ even (resp. odd),
if $\bsb f$, $|\bsb f|\geq 2$, does not have the form $0^\ell$ nor $(q-1)0^\ell$ 
(resp.  the form $0^\ell$) for some $\ell\geq 1$,
then for any $\bsb p\in A_q^*(\bsb f)$, at least one among 
the $\prec$-first and the $\prec$-last word in $\bsb p|A_q^\infty(\bsb f)$
(resp. the $\triangleleft$-first and the $\triangleleft$-last word 
in $\bsb p|A_q^\infty$) has ultimate period $0$.
\end{Rem}

\subsection{Forbidden factor inducing zero periodicity}
\medskip

Here we characterize the first and the last words in $\bsb p|A_q^\infty(\bsb f)$
when the forbidden factor $\bsb f$ induces zero periodicity; the resulting ultimate $0$ periodic words
will be used in the next section.

\begin{Pro}
\label{pro:zerosuffix_end0}
Let $q\geq 2$ be even, $\bsb f\in A_q^+\setminus U_q$ be a forbidden factor 
ending by $0$, $\ell\geq 1$ be the length of the maximal $0$ suffix of $\bsb f$, 
and $\bsb p\in A_q^*(\bsb f)$.
If $\bsb a$ is the $\prec$-first or the $\prec$-last
word in $\bsb p|A_q^\infty(\bsb f)$, then $\bsb a$ has the form
$$
\bsb p\bsb r 0^\infty,
$$
where 
\begin{enumerate}
\item $\bsb r=\epsilon$ or  
      $\bsb r=0^i1(q-1)$ for some $i$, $0\leq i\leq \ell-1$, if $\bsb a$ is 
      the $\prec$-first (resp. $\prec$-last) word in $\bsb p|A_q^\infty(\bsb f)$
      and $\bsb p$ has even (resp. odd) parity, or
\item $\bsb r=q-1$ or $(q-1)0^{\ell-1}1(q-1)$ if $\bsb a$ is 
      the $\prec$-first (resp. $\prec$-last) word in $\bsb p|A_q^\infty(\bsb f)$
      and $\bsb p$ has odd (resp. even) parity.
\end{enumerate}
\end{Pro}
\ccoomm{
\proof We prove the first point, the second one being similar. 
Let $\bsb a$ be the $\prec$-first (resp. $\prec$-last) word in $\bsb p|A_q^\infty(\bsb f)$
with $\bsb p$ having  even (resp. odd) parity. Let also suppose that $\bsb r$
has not the form prescribed in point 1.
Reasoning as in the proof of Proposition \ref{pro:end_0_even} it follows 
that $0^i1(q-1)0^{\ell-1}1(q-1)0^{\ell-1}$ is a prefix of $\bsb r$, for some $i$, 
$0\leq i\leq \ell-1$, and finally, by Lemma \ref{comb_w} that $\bsb f\in U_q$, 
which leads to a contradiction.
}
\endproof

The proof of the next proposition is similar to that of Proposition \ref{pro:zerosuffix_end0}.

\begin{Pro}
\label{zerosuffix_q_odd_end0}
Let $q\geq 3$ be odd, $\bsb f\in A_q^+\setminus V$ be a forbidden factor ending by $0$,
$\ell\geq 1$ be the length of the maximal $0$ suffix of $\bsb f$, 
and $\bsb p\in A_q^*(\bsb f)$.
If $\bsb a$ is the $\triangleleft$-first or the $\triangleleft$-last
word in $\bsb p|A_q^\infty(\bsb f)$, then $\bsb a$ has the form
$$
\bsb p\bsb r 0^\infty,
$$
where 
\begin{enumerate}
\item $\bsb r=\epsilon$ or  
      $\bsb r=0^i1$ for some $i$, $0\leq i\leq \ell-1$, if $\bsb a$ is 
      the $\triangleleft$-first (resp. $\triangleleft$-last) word in $\bsb p|A_q^\infty(\bsb f)$
      and $\bsb p$ has even (resp. odd) parity, or
\item $\bsb r=q-1$ or $(q-1)0^{\ell-1}1$ if $\bsb a$ is 
      the $\triangleleft$-first (resp. $\triangleleft$-last) word in 
      $\bsb p|A_q^\infty(\bsb f)$
      and $\bsb p$ has odd (resp. even) parity.
\end{enumerate}
\end{Pro}

It is routine to check the next two propositions.

\begin{Pro}
\label{oth_pro}
Let $q\geq 3$, $\bsb f\in A_q^+\setminus W_q$ be a forbidden factor ending by $q-1$,
and $\bsb p\in A_q^*(\bsb f)$.
If $\bsb a$ is the first or the last
word in $\bsb p|A_q^\infty(\bsb f)$ with respect to the appropriate order, then 
$\bsb a$ 
has the form
$$
\bsb p\bsb r 0^\infty,
$$
where $\bsb r$ is either $\epsilon$, or $q-1$, or $(q-2)(q-1)$.
\end{Pro} 

As mentioned in Remark \ref{easy}, forbidden factors ending by other symbol
than $0$ or $q-1$ induce zero periodicity, and we have
the following proposition.

\begin{Pro}
\label{pro:zerosuffix_end_other}
If $\bsb f\in A_q^*$ is a  forbidden factor 
that does not end by $0$ nor by $q-1$, then for any ${\bsb p}\in A_q^*({\bsb f})$, 
with respect to the appropriate order, both the first and the last 
word in ${\bsb p}|A_q^\infty({\bsb f})$ have the form:
$$
{\bsb p}{\bsb r}0^\infty,
$$
where $\bsb r$ is either $\epsilon$ or $q-1$.
\end{Pro}

We will see later that Propositions \ref{pro:zerosuffix_end0} to 
\ref{pro:zerosuffix_end_other} above describe sufficient (but not a necessary) 
conditions for the Graycodeness of $A_q^n(\bsb f)$. 

We conclude this section by the next corollary which 
summarizes the results in 
Remark \ref{easy} and Theorems \ref{the:end_zero_even}, \ref{the:end_zero} 
and \ref{the:end_max}, and we will refer it later.

\begin{Co}
\label{Cor}
The forbidden factor $\bsb f\in A_q^*$ induces zero periodicity if and only if either:
\begin{itemize}
\item $\bsb f$ does not end by $0$ nor by $q-1$, or
\item $q=2$ and $\bsb f\not\in U_2$, or
\item $q\geq 4$ is even and ${\bsb f}\not\in U_q\cup W_q$, or
\item $q\geq 3$ is odd and ${\bsb f}\not\in  V\cup W_q$.
\end{itemize}
\end{Co}

\section{Gray codes}
\label{Sec_Gray_codes}

In this section we show that for forbidden factors
$\bsb f$ inducing zero periodicity on $A_q^\infty$ (as stated in Corollary \ref{Cor})
consecutive words---in $\prec$ order for  $q$ even,
or $\triangleleft$ order for $q$ odd---in $A_q^n({\bsb f})$, 
beyond the common prefix, have all symbols $0$, except the first few of them; and
this ensures that the set $A_q^n({\bsb f})$ listed in 
an appropriate order is a Gray code.

Nevertheless, the property of $\bsb f$ to induce zero periodicity
is not a necessary condition.
Indeed, listing the set $A_q^n(\bsb f)$ in $\prec$ order 
with: 
\begin{itemize}
\item $\bsb f=0^\ell$ for any $q$ (not necessarily even), or 
\item $\bsb f=(q-1)0^\ell$ for  $q$ even,
\end{itemize}
where $\ell\geq 1$, yields a $1$-Gray code, despite $\bsb f\in U_q$ (and so, $\bsb f$ does not induce 
zero periodicity for $q$ even).
%
This particular cases are discussed in Section \ref{sec:particular},
and we show that such factors $\bsb f$, $\bsb f\geq 2$, are the only ones giving Gray codes for forbidden factors
not inducing zero periodicity.
In particular, the Gray code obtained for $A_q^n(0^\ell)$ is one of those defined in 
\cite{BBPV_2013} as a generalization of a Gray code in \cite{Vaj_2001}.
Finally, for forbidden factors $\bsb f$ for which $\prec$ nor $\triangleleft$
does not produce Gray codes on $A_q^n({\bsb f})$, we give simple transformations
of $\bsb f$, and eventually obtain Gray codes for $A_q^n({\bsb f})$
(in order other than $\prec$ or $\triangleleft$).

\medskip

We will make use later of the following property of
forbidden factors ending by $0$ or $q-1$: for {\it any} $q\geq 2$,
if $\bsb f$ ends by $0$ or $q-1$, then any two consecutive words in $A_q^n(\bsb f)$, in both
$\prec$ and $\triangleleft$ order, 
differ by $1$ or $-1$
in the leftmost position where they differ.

\begin{Pro}
\label{pm1}
Let $q\geq 2$ and $\bsb f\in A_q^+$ be a forbidden factor ending by $0$ 
or $q-1$, 
and $\bsb a=a_1a_2\ldots a_n$ and  $\bsb b=b_1b_2\ldots b_n$ 
be two words in $A_q^n(\bsb f)$, consecutive with respect to $\prec$ 
or $\triangleleft$ order.
If $k$ is the leftmost position where $\bsb a$ and $\bsb b$ differ,
then $b_k=a_k+1$ or $b_k=a_k-1$.
\end{Pro}
\ccoomm{
\proof
If $\bsb f$ ends by $0$ let us suppose that $b_k<a_k-1$. It follows that $\bsb f$
is a suffix of $a_1a_2\ldots (a_k-1)$, 
so $a_k-1=0$ and thus $b_k<0$, which is a contradiction.
The proof when $b_k>a_k+1$ or when $\bsb f$ ends by $q-1$
is similar.
}
\endproof

\subsection{Factors inducing zero periodicity}
\label{sec:general}

We show that for factors $\bsb f$ as in  Corollary 
\ref{Cor} the set $A_q^n(\bsb f)$ listed in $\prec$ or $\triangleleft$ order
is a Gray code.

\begin{Pro}
\label{pro:general_other}
If $q$ is even (resp. odd) and ${\bsb f}\in A_q^+$ does not end by $0$ nor 
 $q-1$, then $A^n_q(\bsb f)$, $n\geq 1$, listed in $\prec$ (resp. $\triangleleft$) order
is a $2$-adjacent Gray code.
\end{Pro}
\ccoomm{
\proof
Let $\bsb a,\bsb b\in A_q^n(\bsb f)$,
$\bsb a=a_1a_2\ldots a_n$ and $\bsb b=b_1b_2\ldots b_n$
be two consecutive words with respect to the appropriate order,
and $k$ be the leftmost position where they differ.
Since $\bsb f$ does not end by $0$ nor by $q-1$, it follows that $q\geq 3$, and considering the 
definitions of $\prec$ and $\triangleleft$ order,
we have in both cases (see Remark \ref{G_VoRGCO})
$\{a_{k+1},b_{k+1}\}\subset\{0,q-1\}$ and $a_{k+2}\ldots a_n=b_{k+2}\ldots b_n=0^{n-k-1}$.
In any case, $\bsb a$ and $\bsb b$ differ in position $k$ and possibly 
in position $k+1$.
}
\endproof

\begin{The}
\label{the:general_end0_qEven}
If $q\geq 2$ is even, ${\bsb f}\in A_q^+\setminus U_q$ ends by $0$,
and $\ell$ is the length of the maximal $0$ suffix of~$\bsb f$,
then $A^n_q(\bsb f)$, $n\geq 1$, listed in $\prec$ order is an at most 
$(\ell+2)$-close $3$-Gray code.
\end{The}
\ccoomm{
\proof
Let $\bsb a,\bsb b\in A_q^n(\bsb f)$,
$\bsb a=a_1a_2\ldots a_n$ and $\bsb b=b_1b_2\ldots b_n$
be two consecutive words with respect to the appropriate order, and $k$ be the leftmost 
position where $\bsb a$ and $\bsb b$ differ.
By Proposition \ref{pm1}, $b_k=a_k+1$ or $b_k=a_k-1$
and so the prefixes $\bsb a'=a_1a_2\ldots a_k$ and $\bsb b'=b_1b_2\ldots b_k$ have 
different parity. Two cases arise according to the parity of $\bsb a'$.

\noindent
$\bullet$ $\bsb a'$ has even parity, and so $\bsb b'$
has odd parity. By point 2 of Proposition \ref{pro:zerosuffix_end0}
$$\bsb a=\bsb a'\bsb x,$$ and 
$$\bsb b=\bsb b'\bsb y,$$
with $\bsb x$ and $\bsb y$ being the $n-k$ prefixes of $\bsb r0^\infty$ and of
$\bsb r'0^\infty$, where $\{\bsb r, \bsb r'\}\subset\{(q-1),(q-1)0^{\ell-1}1(q-1) \}$. 
Thus $\bsb a$ and $\bsb b$ differ in position $k$ and possibly in 
positions $k+\ell+1$ and $k+\ell+2$ if $k+\ell+1\geq n$.

\noindent
$\bullet$ $\bsb a'$ has odd parity, and so $\bsb b'$
has even parity. By point 1 of Proposition \ref{pro:zerosuffix_end0} either
\begin{itemize}
\item[(i)] $a_{k+1}a_{k+2}\ldots a_n=b_{k+1}b_{k+2}\ldots b_n=0^{n-k}$, or
\item[(ii)] at least one of $a_{k+1}a_{k+2}\ldots a_n$ or $b_{k+1}b_{k+2}\ldots b_n$
        is the $n-k$ prefix of a word of the form  $0^i1(q-1)0^\infty$.
\end{itemize}
In case (i) $\bsb a$ and $\bsb b$ differ only in position $k$. And in case 
(ii) let us suppose that  $a_{k+1}a_{k+2}\ldots a_n$ is the length $n-k$ prefix of 
$0^i1(q-1)0^\infty$
(the corresponding case for $b_{k+1}b_{k+2}\ldots b_n$ being similar).
Considering that $b_k=a_k+1$ or $b_k=a_k-1$
it follows that $b_{k+1}b_{k+2}\ldots b_n$ is the length $n-k$ prefix of $0^\infty$ and 
so $\bsb a$ and $\bsb b$ differ in positions $k$, and 
(possibly) $k+i+1$ and $k+i+2$.

In any case, $\bsb a$ and $\bsb b$ differ in at most three positions
which are at most $\ell+2$ apart from each other.
}
\endproof

\begin{Exam}
The words $00230130$ and $00330000$ are consecutive in $A_4^8(2300)$ listed in 
$\prec$ order. They differ in $3$ positions which are $4$-close, and are 
in the worst case since the list is a $4$-close $3$-Gray code.
\end{Exam}

Considering the possible values of $\bsb r$ in Proposition \ref{zerosuffix_q_odd_end0}
it is easy to see that for $\bsb f\notin V$ ending by $0$ and 
$q$ odd, the set $A^n_q(\bsb f)$ listed in $\triangleleft$ order is a $4$-Gray code.
The next theorem gives a more restrictive result.

\begin{The}
\label{the:general_end0_qOdd}
If $q\geq 3$ is odd, ${\bsb f}\in A_q^+\setminus V$ ends by $0$, and $\ell$ is the 
length of the maximal $0$ suffix of~$\bsb f$,
then $A^n_q(\bsb f)$, $n\geq 1$, listed in $\triangleleft$ order is an at most 
$(\ell+1)$-close $3$-Gray code.
\end{The}
\ccoomm{
\proof
Let $\bsb a,\bsb b\in A_q^n(\bsb f)$,
$\bsb a=a_1a_2\ldots a_n$ and $\bsb b=b_1b_2\ldots b_n$
be two consecutive words, in $\triangleleft$ order, and $k$ be the leftmost 
position where $\bsb a$ and $\bsb b$ differ. If $\bsb a'$ and $\bsb b'$ are the length $k$ 
prefix of $\bsb a$ and $\bsb b$, by Proposition \ref{zerosuffix_q_odd_end0}
$$\bsb a=\bsb a'\bsb x,$$ and 
$$\bsb b=\bsb b'\bsb y,$$
with $\bsb x$ and $\bsb y$ being the  $n-k$ prefixes of $\bsb r0^\infty$ and of
$\bsb r'0^\infty$, where 
$\{\bsb r, \bsb r'\}\subset
\{\epsilon, 0^i1,(q-1),(q-1)0^{\ell-1}1\}$, for some $i$, $0\leq  i\leq \ell-1$.
The statement holds by showing that
$\{\bsb r, \bsb r'\}\subset\{0^i1,(q-1)0^{\ell-1}1\}$
is not possible. Indeed, let us suppose that 
$\bsb r=0^i1$ for some $i$, $0\leq  i\leq \ell-1$,
and $\bsb r'=(q-1)0^{\ell-1}1$
(the case $\bsb r=(q-1)0^{\ell-1}1$ and $\bsb r'=0^i1$ being similar). This happens when both 
$\bsb a'$ and $\bsb b'$ have both odd parity. 
By Proposition \ref{pm1}, $b_k=a_k+1$ or $b_k=a_k-1$, and 
since $a_1a_2\ldots a_{k-1}=b_1b_2\ldots b_{k-1}$
it follows that $a_k=1$ and $b_k=0$.
Since $\bsb r'=(q-1)0^{\ell-1}1$, the factor $\bsb f$ must end 
by $(q-1)0^{\ell}$ and since $a_k=1$ it follows that $\bsb r=\epsilon$,
which leads to a contradiction.
}
\endproof

\begin{Exam}
By Theorem \ref{the:general_end0_qOdd}, the sets  $A_5^9(31000)$ and $A_5^9(24000)$ listed in
$\triangleleft$ order are 4-close $3$-Gray codes. However, it is easy to check that in particular,
$A_5^9(31000)$ is a  $3$-close $3$-Gray code, and $A_5^9(24000)$ is $4$-close $2$-Gray code.
For example:
\begin{itemize} 
\item the words $001304000$ and $001310010$ are consecutive in $A_5^9(31000)$  when listed in 
$\triangleleft$ order; 
they differ in $3$ positions which are $3$-close; and
\item the words $001140000$ and $001240010$ are consecutive in $A_5^9(24000)$  when listed in 
$\triangleleft$ order; 
they differ in $2$ positions which are $4$-close. 
\end{itemize}

\end{Exam}

\begin{The}
\label{the:general_max}
If $q$ is even (resp. odd) and ${\bsb f}\in A_q^+
\setminus W_q$ ends by $q-1$,
then $A^n_q(\bsb f)$ listed in $\prec$ order (resp. $\triangleleft$ order) is a 
$2$-close $3$-Gray code (that is, a $3$-adjacent Gray code).
\end{The}
\ccoomm{
\proof
Let $k$ be the leftmost position where two consecutive words
$\bsb a=a_1a_2\ldots a_n$ and  $\bsb b=b_1b_2\ldots b_n$,
in $A^n_q(\bsb f)$ differ. Exhausting the possible values of $a_k$ and $b_k$, 
and since ${\bsb f}\not\in W_q$ ends by $q-1$ it follows that $a_i=b_i=0$ for all $i>k+2$
(see also Proposition \ref{oth_pro}).
}
\endproof

\subsection{Particular cases}
\label{sec:particular}

As mentioned before, there are two cases when $\bsb f\in U_q$, $q\geq2$ and even,
but $A_q^n(\bsb f)$ listed in $\prec$ order is a Gray code; 
these are $\bsb f=0^{\ell}$ and $\bsb f=(q-1)0^\ell$, $\ell\geq 1$. 
Moreover, it turns out that $A_q^n(0^\ell)$, $q\geq3$ and odd, also gives Gray code 
if listed in $\prec$ order. 
Similar phenomenon does not occur for $\bsb f\in V$, i.e.,
the set $A_q^n(\bsb f)$ listed in $\triangleleft$ order is not a Gray code for any 
$\bsb f\in V$, $|\bsb f|\geq 2$ and $q\geq3$ odd, see for example 
Remark \ref{except_q_odd}.

\medskip

Before discussing these particular forbidden factors we introduce some notations.

Let $q\ge2$, $\ell\geq 1$, and let define the infinite words:

\begin{equation}
\label{u_v}
\begin{array}{l}
\bsb u = (0^{\ell-1}1(q-1))^\infty,\\
\bsb v = ((q-1)0^{\ell-1}1)^\infty.
\end{array}
\end{equation}

Notice that $\bsb u$ and $\bsb v$ are suffixes 
to each other, and they are related with the infinite words occurring in Proposition 
\ref{pro:end_0_even}. It is easy to see that $\bsb u$ and $\bsb v$ are, respectively,
the $\prec$-first and $\prec$-last word in $A_q^\infty(0^\ell)$ for even $q$; and thus
the length $n$ prefix of $\bsb u$ and $\bsb v$ are, respectively,
the $\prec$-first and $\prec$-last word in $A_q^n(0^\ell)$.

Moreover,  for any $\bsb p\in A_q^k(0^\ell)$ with  $1\leq k\le n$ and $q$ even
\begin{itemize}
\item the $\prec$-first (resp. $\prec$-last) word in
      $\bsb p|A_q^n(0^\ell)$ is $\bsb p\bsb v'$ if $\bsb p$ has an odd
      (resp. even) parity, where $\bsb v'$ is the length $n-k$
      prefix of  $\bsb v$;
\item if $\bsb p$ does not end by $0$, then the $\prec$-first (resp. $\prec$-last) word in
      $\bsb p|A_q^n(0^\ell)$ is $\bsb p\bsb u'$ if $\bsb p$ has an even
      (resp. odd) parity, where $\bsb u'$ is the length $n-k$
      prefix of  $\bsb u$.
\end{itemize}

Now let $q\geq 3$ and odd, $\ell\geq 1$, and let define the infinite words:

\begin{equation}
\label{s_t}
\begin{array}{l}
\bsb s = 0^{\ell-1}1(q-1)^\infty,\\
\bsb t = (q-1)^\infty,
\end{array}
\end{equation}
and $\bsb s$ and $\bsb t$ have similar property as $\bsb u$ and $\bsb v$ 
for $q$ odd and with same $\prec$ order.

\subsubsection*{The case $\bsb f=0^\ell$}
\medskip

\begin{Pro}
\label{all_q_0_k}
For $q\geq 2$ (even or odd), and $\ell,n\geq 1$, the set 
$A_q^n(0^\ell)$ listed in $\prec$ order is a Gray code
where two consecutive words differ in one position and by $1$ or $-1$
in this position.
\end{Pro}
\ccoomm{
\proof
Let $\bsb a$ and $\bsb b$ be two consecutive words, 
in $\prec$ order, in  $A_q^n(0^\ell)$,
$\bsb a'=a_1a_2\ldots a_k$ and $\bsb b'=b_1b_2\ldots b_k$ be the length $k$ prefix of 
$\bsb a$ and $\bsb b$, with $k$ the leftmost position where 
$\bsb a$ and $\bsb b$ differ.

When $q$ is even, with $\bsb u'$ and $\bsb v'$ the length $(n-k)$ prefix of 
$\bsb u$ and $\bsb v$ defined in relation~(\ref{u_v}), we have
\begin{itemize}
\item $\bsb a=\bsb a'\bsb v'$ and  $\bsb b=\bsb b' \bsb v'$ if $\bsb a'$ has an even parity 
      (and so, by Proposition \ref{pm1},
     $\bsb b'$ has odd parity);
\item $\bsb a=\bsb a'\bsb u'$ and $\bsb b=\bsb b'\bsb u'$, elsewhere, since
      $a_k\neq 0$ and $b_k\neq 0$ by 
      considering the parity of the common length $k-1$ prefix of $\bsb a$ and $\bsb b$.
\end{itemize}

Similarly, when $q$ is odd, with $\bsb s'$ and $\bsb t'$ the length $(n-k)$ prefix of 
$\bsb s$ and $\bsb t$ defined in relation~(\ref{s_t}), we have
\begin{itemize}
\item $\bsb a=\bsb a'\bsb t'$ and  $\bsb b=\bsb b'\bsb t'$ 
if $\bsb a'$ has an even parity (given by $\sum_{i=1}^k a_i$);
\item $\bsb a=\bsb a'\bsb s'$ and  $\bsb b=\bsb b'\bsb s'$ 
(since, $a_k\neq 0$ and $b_k\neq 0$ ), elsewhere.
\end{itemize}
In both cases $\bsb a$ and $\bsb b$ differ only in position $k$.
}
\endproof

\subsubsection*{The case $\bsb f=(q-1)0^{\ell}$ for $q$ even}
\medskip
Let $q\geq2$ be even, $1\leq k\le n$, and $\bsb v'$ be the $n-k$ prefix of 
$\bsb v$ defined in relation (\ref{u_v}).
For any $\bsb p\in A_q^k((q-1)0^{\ell})$ with $\ell\geq 1$ and $1\leq k\le n$
\begin{itemize}
\item the $\prec$-first (resp. $\prec$-last) word in
      $\bsb p|A_q^n((q-1)0^{\ell})$ is $\bsb p\bsb v'$ if $\bsb p$ has odd  
      (resp. even) parity;
\item if $\bsb p$ does not end by $0$ nor by $q-1$, then the $\prec$-first 
      (resp. $\prec$-last) word in $\bsb p|A_q^n((q-1)0^{\ell})$ is $\bsb p0^{n-k}$ 
      if $\bsb p$ has even (resp. odd) parity.
\end{itemize}

\begin{Pro}
\label{even_q_max}
For $q\geq 2$ even, and $\ell,n\geq 1$, the set 
$A_q^n((q-1)0^\ell)$ listed in $\prec$ order is a Gray code
where two consecutive words differ in one position and by $1$ or $-1$
in this position.
\end{Pro}
\ccoomm{
\proof 
Let $\bsb a$ and $\bsb b$ be two consecutive words, in $\prec$ order,
in $A_q^n((q-1)0^{\ell})$, $\bsb a'=a_1a_2\ldots a_k$ and $\bsb b'=b_1b_2\ldots b_k$ be the length $k$ prefix of 
$\bsb a$ and $\bsb b$ with $k$ the leftmost position where 
$\bsb a$ and $\bsb b$ differ. 

If $\bsb a'$ has even parity
(and so, by Proposition \ref{pm1}, $\bsb b'$ has odd parity), then
by the above considerations $\bsb a=\bsb a'\bsb v'$ and  
$\bsb b=\bsb b'\bsb v'$.

If $\bsb a'$ has odd parity, by considering the parity
of the common length $k-1$ prefix of $\bsb a$ and $\bsb b$
it follows that $a_k\neq q-1$ and $b_k\neq q-1$, and again, by 
the above considerations we have $\bsb a=\bsb a'0^{n-k}$ and  
$\bsb b=\bsb b'0^{n-k}$.

In both cases $\bsb a$ and $\bsb b$ differ only in position $k$.
}
\endproof

\subsection{Factors preventing Graycodeness}

A consequence of the next remark and proposition, 
is Corollary \ref{Co_other} below. Proposition \ref{ult_perod} sounds like Remark \ref{useless}
and says that if $\bsb f$, $|\bsb f|\geq 2$
($|\bsb f|=1$ being trivial), does not induce zero periodicity (see Corollary \ref{Cor}), and 
it is not in one of the two particular cases above, then consecutive words, with 
respect to the appropriate order, in $A_q^n(\bsb f)$ can differ in arbitrary many positions
for enough large $n$. 
One of these particular cases is explained below.

\begin{Rem}
\label{except_q_odd}
For $q\geq 3$ and odd, $\ell\geq 2$ and $\bsb f=0^{\ell}$, the set $A_q^n(\bsb f)$ listed in $\triangleleft$-order 
is not a Gray code. Indeed, for example, the words $02\bsb z'$ and $1\bsb z''$ are consecutive 
in $\triangleleft$-order in $A_q^n(\bsb f)$, where $\bsb z'$ and  $\bsb z''$
are appropriate length prefixes of $(0^{\ell -1}1)^\infty$, and they differ in arbitrary many positions
for enough large $n$.
\end{Rem}



\begin{Pro}
\label{ult_perod}
Let $\bsb f\in A_q^+$, $q\geq 2$ and $|\bsb f|\geq 2$, be a forbidden factor
not inducing zero periodicity, other than $0^\ell$ or $(q-1)0^\ell$, 
$\ell\geq 1$.
Let also $\bsb a$ and $\bsb b$ be two consecutive words, in appropriate order, in 
$A_q^n(\bsb f)$, $n\geq 1$, and $k$ the leftmost position where $\bsb a$ and $\bsb b$ differ.
If 
\begin{itemize}
\item $\bsb a'=a_1a_2\ldots a_k$ and $\bsb b'=b_1b_2\ldots b_k$ are, respectively, 
      the length $k$ prefix of $\bsb a$ and $\bsb b$,
      and
\item $\bsb a''$ and $\bsb b''$ are, respectively, 
      the last word in $\bsb a'|A_q^\infty(\bsb f)$ 
      and the first word in $\bsb b'|A_q^\infty(\bsb f)$, in appropriate order,     
\end{itemize}
then at most one among $\bsb a''$ and $\bsb b''$
does not have ultimate period $0$.
\end{Pro}
\ccoomm{
\proof
Since $\bsb f$ does not induce zero periodicity, we 
prove the statement according to $\bsb f$ belongs to $U_q$, $V$ or $W_q$
(see Corollary \ref{Cor}),
and supposing that $\bsb a''$ does not have ultimate period $0$ 
(the corresponding case for $\bsb b''$ being similar). 

\noindent
If $\bsb f\in U_q$, $q\geq 2$ and $\bsb f$ does not have the form $0^\ell$ nor $(q-1)0^\ell$:

\begin{itemize}
\item When $\bsb a'$ has odd parity, since $a_k$ must be a symbol
      of $\bsb f$, it follows that $a_k\in\{0,1,q-1\}$.
       From the parity of $\bsb a'$, it follows that 
       $a_k=0$ implies that $b_k=a_k-1$, and  $a_k=q-1$ that $b_k=a_k+1$, which are not possible,
       and necessarily $a_k=1$. Thus, either $\bsb a''=\bsb a'0^\infty$ (which is a contradiction
       with the non-zero periodicity of $\bsb a$) or $\bsb b''=\bsb b'0^\infty$.
\item  When $\bsb a'$ has even parity, then 
       $\bsb a''=\bsb a'\bsb v$ and since $b_k\neq a_k$,
       $\bsb b''=\bsb b'(q-1)0^\infty$, with $\bsb v$ defined in relation 
       (\ref{u_v}).
\end{itemize}

\noindent
If $\bsb f\in V$, $q\geq 3$ and odd, and $\bsb f$ does not have the form $0^\ell$:

\begin{itemize}
\item $\bsb a'$ can not have even parity, otherwise $\bsb a''=\bsb a'(q-1)0^\infty$,
      which is a contradiction with the non-zero periodicity of $\bsb a$;
\item When $\bsb a'$ has odd parity, the symbol $a_k$ must be one of the forbidden factor, so 
      $a_k\in \{0,1\}$. But $a_k=0$, implies $b_k=a_k-1$,
      which again is not possible; and $a_k=1$ implies $b_k=0$, and so $\bsb b''=(q-1)0^\infty$,
      which does not contain the factor $\bsb f$ if it is different from $0^\ell$.
\end{itemize}

\noindent
Finally, when $\bsb f\in W_q$, $q\geq 3$ (even or odd) then
$\bsb a''=\bsb a'(q-2)^\infty$ and  $\bsb b''$ is either $\bsb b'0^\infty$ 
(this can occur if $q$ is odd) or $\bsb b'(q-1)0^\infty$.
}
\endproof

Table \ref{summary} summarizes the cases occurring in Proposition \ref{ult_perod}.

A consequence of Remark \ref{except_q_odd}, Propositions \ref{all_q_0_k} to \ref{ult_perod} 
and Corollary \ref{Cor}, is the corollary below.

\begin{Co}
\label{Co_other}
$ $
\begin{itemize}
\item For even $q\geq 2$ and $|\bsb f|\geq 2$, the set $A_q^n(\bsb f)$ listed in 
      $\prec$ order is a 
      Gray code for any $n\geq1$ if and only if 
      $\bsb f\in \{0^\ell,(q-1)0^\ell\}_{\ell\geq1}\cup W_2\cup (A_q^*\setminus  (U_q\cup
      W_q))$.
\item For odd $q\geq 3$ and $|\bsb f|\geq 2$,  the set $A_q^n(\bsb f)$ listed in 
       $\triangleleft$ order is a Gray code for any $n\geq1$ 
      if and only if $\bsb f\in  A_q^*\setminus (V\cup W_q)$.
      
\end{itemize}
\end{Co}

\begin{table}[h]
\centering
\begin{tabular}{| >{\centering\arraybackslash}m{1.5cm}|>{\centering\arraybackslash}m{1cm}|>{\centering\arraybackslash}m{4cm} |>{\centering\arraybackslash}m{3.cm}|>{\centering\arraybackslash}m{3cm}|}
\hline
 $q$ & Order relation &The set for forbidden factor $\bsb f$ & The set of ultimate periods of the last 
 word in $\bsb a|A_q^\infty(\bsb f)$ and  the first one in $\bsb b|A_q^\infty(\bsb f)$ & Graycodeness  of $A_q^n(\bsb f)$\\ 
 \hline\hline
  even & $\prec$ & $U_q\setminus\{0^{\ell},(q-1)0^\ell\}_{\ell\geq1}$& $\{1(q-1)0^{\ell-1}, 0\}$ 
 & Not Gray code\\[.5ex]
\hline 
even & $\prec$ &  $\{0^{\ell},(q-1)0^\ell\}_{\ell\geq1}$& $\{1(q-1)0^{\ell-1}\}$
& 1-Gray code\\[.5ex]
\hline 
odd & $\prec$ & $\{0^\ell\}_{\ell\geq1} $ & $\{(q-1)\}$
& 1-Gray code\\[.5ex]
 \hline \hline
odd & $\triangleleft$ & $V\setminus\{0^\ell\}_{\ell\geq1}$ &$\{10^{\ell-1}, 0\}$ 
& Not Gray code \\[.5ex]
\hline 
odd & $\triangleleft$ & $\{0^\ell\}_{\ell\geq2} $ & $\{10^{\ell-1}\}$
& Not Gray code\\
\hline\hline
  {$q\geq3$} even (resp. odd)& $\prec$ (resp. $\triangleleft$) & $W_q\cap A_q^{\geq 2}$ & $\{(q-2), 0\}$ & Not Gray code\\
\hline
 
\end{tabular}
\caption{
\label{summary}
The Graycodeness  of $A_q^n(\bsb f)$ listed in appropriate order
together with the  ultimate periods of the last word in $\bsb a|A_q^\infty(\bsb f)$
and the first word in $\bsb b|A_q^\infty(\bsb f)$, when at least one of them
does not have ultimate period $0$, and $\bsb a$ and $\bsb b$ are consecutive 
words; and $A_q^{\geq 2}$ is the set of words on $A_q$ of length at least two.
These summarize Propositions \ref{all_q_0_k} to \ref{ult_perod}, and 
Corollary \ref{Co_other}.}
\end{table}

\subsection{Obtaining Gray code if $\bsb f$ does not induce zero periodicity and 
beyond the particular cases}

According to the previous results, if the forbidden factor $\bsb f$
does not induce zero periodicity, then the set $A_q^n(\bsb f)$ listed
in $\prec$ or $\triangleleft$ order is not a Gray code, except for the two 
particular cases in Section \ref{sec:particular}.
Now we show how a simple transformation allows to define Gray codes,
with the same Hamming distance and 
closeness properties as for factors that induce zero periodicity,
when $\bsb f$ does not have this property.
By Theorems \ref{the:end_zero_even}, \ref{the:end_zero} and \ref{the:end_max}, the last symbol of 
a factor that does not induce zero periodicity is either $0$, or $q-1$ when $q\geq 3$. 

Let define the transformation $\phi:A_q\rightarrow A_q$ depending on $\bsb f$ as 
\begin{itemize}
\item when the last symbol of $\bsb f$ is $0$, then $\phi(0)=1$, $\phi(1)=0$,
      and $\phi(x)=x$ if $x\not\in\{0,1\}$; and 
\item when the last symbol of $\bsb f$ is $q-1$ (and so, $q\geq 3$), then $\phi(q-2)=q-1$, 
      $\phi(q-1)=q-2$, and $\phi(x)=x$ if $x\not\in\{q-2,q-1\}$.
\end{itemize}

In both cases, $\phi$ is an involution, that is, $\phi^{-1}=\phi$.
By abuse of notation, 
for ${\bsb w}\in A^*_q$, $\phi({\bsb w})$ is the word obtained from ${\bsb w}$
by replacing each of its symbols $x$ by $\phi(x)$, and for a list ${\mathcal L}$ of words,
$\phi({\mathcal L})$ is the list obtained from ${\mathcal L}$ by replacing each word ${\bsb w}$
in ${\mathcal L}$ by $\phi({\bsb w})$.

If ${\bsb f}$ is a forbidden factor that does not induce zero periodicity, 
then $\phi({\bsb f})$ does not end by $0$ nor by $q-1$, and so it 
induces zero periodicity, see Remark \ref{easy}.
In this case $\phi({\mathcal L})$ is a Gray code for the set $A_q^n({\bsb f})$, 
where ${\mathcal L}$ is the set $A_q^n(\phi({\bsb f}))$ listed in $\prec$ order
for $q$ even, and in $\triangleleft$ order for $q$ odd.

\section{Algorithm considerations}
\label{sec:algorithm}

Here we give a generating algorithm  for the set $A_q^n(\bsb f)$, $n\geq 1$, 
for any forbidden factor $\bsb f\in A_q^\ell$, $\ell\geq 2$ (the case $\ell=1$ being trivial). 
This generating algorithm produces recursively prefixes of words in  
$A_q^n(\bsb f)$, in $\prec$ order if $q$ is even, or in $\triangleleft$ order 
if $q$ is odd, and in particular, it generates 
the previously discussed Gray codes for $A_q^n(\bsb f)$. 
We will show that this algorithm is efficient, except for the trivial 
factors of the form $00\ldots 01$ or $11\ldots 10$, for which a simple
transformation of them makes the generating algorithm efficient.

The generating procedure {\sc GenAvoid} in Figure \ref{fig:algorithm_gen} expands recursively
a current generated prefix $w_1w_2\ldots w_{k-1}$
($k$ being the first parameter of {\sc GenAvoid}) to 
$w_1w_2\ldots w_{k-1}j$, with $j$ covering the alphabet $A_q$ in increasing 
or decreasing order, according to the value of $dir\in \{0,1\}$, 
the second parameter of the procedure, which is the parity of the word 
$w_1w_2\ldots w_{k-1}$. 
Moreover, when the length $(\ell-1)$ prefix of $\bsb f=f_1f_2\ldots f_\ell$
is a suffix of $w_1w_2\ldots w_{k-1}$, 
the value $f_\ell$ is skipped for $j$ in order not to produce the 
forbidden factor. 
To do this efficiently, the third parameter, $i$, of procedure {\sc GenAvoid} 
is the length of the maximal prefix of the  forbidden factor $\bsb f$
which is also a  suffix of the current generated word $w_1w_2\ldots w_{k-1}$;
and $h$ in this procedure is the length of the maximal 
suffix of $w_1w_2\ldots w_{k-1}j$ which is also a prefix of $\bsb f$, and 
it is given by $M_{i,j}$. 
So, when $h$ is equal to $\ell$ (the length of the forbidden factor),
the current value of $j$ is skipped for the prefix expansion.

Now we explain in more details the array $M$ 
used by algorithm  {\sc GenAvoid}.
For a forbidden factor $\bsb f=f_1f_2\ldots f_\ell\in A_q^\ell$, 
the $\ell\cdot q$ size two dimensional array $M$ is defined as:
for $i\in \{0,1,\ldots , \ell-1\}$ and $j\in\{0,1,\ldots ,q-1\}=A_q$,
$M_{i,j}$ is the length of the maximal suffix of $f_1f_2\ldots f_ij$
which is also a prefix $\bsb f$.
For instance, for $q=4$ and $\bsb f=012011\in A_4^6$
we have

$$
M=
\left[
\begin{array}{cccc}
 1 & 0 & 0 & 0  \\
 1 & 2 & 0 & 0  \\
 \bsb 1 &  \bsb 0 &  \bsb 3 & 0  \\
 4 & 0 & 0 & 0  \\
 1 & 5 & 0 & 0  \\
 1 & 6 & 3 & 0 
\end{array}
\right],
$$
and, for example (see the entries in boldface in $M$)
\begin{itemize}
\item $M_{2,0}={\bsb 1}$, since the length of the longest suffix of $f_1f_20=010$ which is a
      prefix of $\bsb f$ is 1,
\item $M_{2,1}={\bsb 0}$, since there is no suffix of $f_1f_21=011$ which is a 
      prefix of $\bsb f$,
\item $M_{2,2}={\bsb 3}$, since $f_1f_22=012$ (of length 3) is a prefix of $\bsb f$.
\end{itemize}

\bigskip

\begin{figure}[!h]

\begin{center}
\begin{tabular}{|c|}
\hline
\begin{minipage}[c]{.46\linewidth}
\begin{tabbing}\hspace{0.9cm}\=\hspace{0.7cm}\=\hspace{1.4cm}\=\hspace{2.9cm}
            \=\hspace{2.9cm}\kill
{\bf procedure} {\sc GenAvoid}($k$, $dir$, $i$)\\
 {\bf if} $k=n+1$ {\bf then} type;\\
  {\bf else} \> {\bf if}  $dir=0$ {\bf then} $\mathcal S:=\langle0,\ldots,q-1\rangle$; 
	        {\bf else} $\mathcal S:=\langle q-1,\ldots,0\rangle$;\\
             \> {\bf for} $j$ in $\mathcal S$ \\
             \>        \>$h:=M[i,j]$;\\        
             \>        \>{\bf if} $h\neq \ell$ \> {\bf then}\\
             \>        \>                      \> $w[k]:=j$; $m:=(dir+j)\mod 2$;\\	         
             \>	       \>                      \> {\bf if} $q$ is odd {\bf and} $j\ne 0$ {\bf then} $m:=(m+1)\mod 2$;\\
             \>        \>                      \> {\sc GenAvoid}($k+1$, $m$, $h$);
\end{tabbing}
\end{minipage}
\\
\hline
\end{tabular}
\caption{
Algorithm producing the set  $A_q^n(\bsb f)$, listed in $\prec$ order if $q$ is 
even or in $\triangleleft$ order if $q$ is odd. 
The initial call is {\sc GenAvoid$(1,0,0)$},  
and it uses array $M$, initialized in a preprocessing step 
by {\sc MakeArray}; and $\mathcal S$ is the list of symbols in the 
alphabet $A_q$ in increasing or decreasing order.
\label{fig:algorithm_gen}}
\end{center}
\end{figure}

The array $M$ is initialized, in an $O(\ell\cdot q)$ time preprocessing  step, by 
procedure {\sc MakeArray} in Figure \ref{fig:algorithm_make_array},
which in turn uses array $b=b_0b_1b_2\ldots b_\ell$, the {\it border} array of $\bsb f$
defined as (see for instance \cite{Lot_2004}): $b_i$, $0\leq i\leq \ell$, is the length of the border of 
$f_1f_2\ldots f_i$,  that is, the length of the longest factor
which is both a proper prefix and a proper suffix of $f_1f_2\ldots f_i$;
and by convenience $b_0=-1$. For example if $\ell=8$ and $\bsb f=01001010$, then  
$b_0b_1\ldots b_8=-100112323$;
and for instance, $b_5=2$ since $01$ is the longest proper prefix 
which is also a suffix of $f_1f_2\ldots f_5=01001$.
Actually, the border array $b$ is a main ingredient for
Knuth-Morris-Pratt word matching algorithm in \cite{Knu_Mor_Pra_1977}
and it is initialized by an $O(\ell)$ time complexity preprocessing step
by procedure {\sc MakeBorder} in Figure \ref{fig:algorithm_make_border}, see again \cite{Lot_2004}.

\begin{figure}[!h]
\begin{center}

\begin{tabular}{|c|}
\hline
\begin{minipage}[c]{.46\linewidth}
\begin{tabbing}
{\bf procedure} {\sc MakeArray$()$}\\
{\bf for} \= $j:=0$ {\bf to} $q-1$\\
	   \>{\bf for} \= $i:=0$ {\bf to} $\ell-1$\\
	   \>     \>{\bf if} $f[i+1]=j$ {\bf then} $M[i,j]:=i+1$;\\
	   \>	\>{\bf else} \= {\bf if} $i>0$ {\bf then} $M[i,j]:=M[b[i],j]$;\\ 
	   \>      \>       \> {\bf else} $M[i,j]:=0$;
\end{tabbing}
\end{minipage}
\\
\hline

\end{tabular}

\caption{Algorithm initializing the array $M$.}
\label{fig:algorithm_make_array}
\end{center}
\end{figure}

\begin{figure}[!h]
\begin{center}
\begin{tabular}{|c|}
\hline
\begin{minipage}[c]{.46\linewidth}
\begin{tabbing}
{\bf procedure} {\sc MakeBorder$()$}\\
   $b[0]:=-1$;\\
   $i:=0$;\\
   {\bf for} \=$j:=1$ {\bf to} $(\ell-1)$ \\
 \>  	$b[j]:=i$;\\
 \>      {\bf whi}\={\bf le} ($i\ge0$ {\bf and} $f[j+1]\ne f[i+1]$)\\
 \> \>	 $i:=b[i]$;\\
  \>      $i:=i+1$;\\
   $b[\ell]:=i$;
\end{tabbing}
\end{minipage}
\\
\hline
\end{tabular}
\caption{Procedure computing the border array $b$ of the 
length $\ell$ forbidden factor $\bsb f$, and used by {\sc MakeArray}.}
\label{fig:algorithm_make_border}
\end{center}
\end{figure}

Before analyzing the time complexity of the generating algorithm {\sc GenAvoid}
we show that, if in the underlying tree induced by recursive calls of {\sc GenAvoid} there are degree-one
successive calls, then $q=2$ and the forbidden factor has the form
$00\ldots01$ or $11\ldots 10$.
See Figure \ref{Fig} for words in $A^n_2(001)$ produced by degree-one consecutive
calls of {\sc GenAvoid}.

\begin{figure}
\begin{center}
\psset{levelsep=17mm,treesep=3mm,radius=-0.5mm} 
\pstree[nodesep=1mm,treemode=R]{\TR{$\cdots$}}
{\pstree{\TR{ $110$}}
{    \pstree{\TR{ $\bsb{1100}$}}
	    {  \pstree{\TR{ $\bsb{11000}$}}
	         {\pstree{\TR{ $\bsb{110000}$}}		   
		   {\pstree{\TR{ $\cdots$}}{}
		   }
		 }
	    }	    
       \pstree{\TR{ $1101$}}
       {\pstree{\TR{ $\cdots$}}{}
        \pstree{\TR{ $\cdots$}}{}
       }  	
}
}
\end{center}
\caption{\label{Fig}
In boldface a `branch' of words produced by consecutive 
degree-one calls in the generating tree of $A^n_2(001)$.
}    
\end{figure}

For a length $\ell\geq 2$ forbidden factor $\bsb f$ 
let $\bsb w\in A_q^*(\bsb f)$ and $i,j\in A_q$ such that $\bsb wij\in A_q^*(\bsb f)$ and:
\begin{itemize}
\item $\bsb wk$ ends by $\bsb f$ for any $k\in A_q$, $k\neq i$, and
\item $\bsb wik$ ends by $\bsb f$ for any $k\in A_q$, $k\neq j$.
\end{itemize}

In other words, when the current word is $\bsb w$ as above, then the call of
{\sc GenAvoid} is a degree-one call (producing $\bsb wi$) which in turn
produces a degree-one call (producing $\bsb wij$).
By the two conditions above, it follows that $q=2$ and $i=j$.
When $i=j=0$, the length $\ell-1$ suffix of $\bsb w$ is equal to the 
$\ell-1$ suffix of $\bsb w0$, which in this case must be $0^{\ell-1}$, and finally
$\bsb f=0^{\ell-1}1$. Similarly, when $i=j=1$, it follows that $\bsb f=1^{\ell-1}0$.

Let now $\bsb f$ be a length $\ell\geq 2$ forbidden factor, and either
$q\geq 3$ or $q=2$ and $\bsb f$ is not $0^{\ell-1}1$ nor $1^{\ell-1}0$.
In this case, by the previous considerations, 
each recursive call of {\sc GenAvoid} is either:
\begin{itemize}
\item a terminal call, or
\item a call producing at least two recursive calls, or
\item a call producing one recursive call, which in turn is in one of the 
two cases above,
\end{itemize}
and by Ruskey's CAT principle in \cite{Ruskey}, it follows that, 
with the previous restrictions on $q$ and $\bsb f$,
{\sc GenAvoid} runs in
constant amortized time, and so is an efficient generating algorithm.

Nevertheless, for the particular factors above,
when $\ell=2$, $A_2^n(1^{\ell-1} 0)$ is trivially the set $\{0^n,0^{n-1}1,0^{n-2}11,\ldots,1^n\}$,
and $A_2^n(0^{\ell-1} 1)$ the set $\{0^n,10^{n-1},110^{n-2},\ldots,1^n\}$.
And for $\ell\geq 3$, both sets  $A_2^n(1^{\ell-1} 0)$ and $A_2^n(0^{\ell-1} 1)$ can be generated 
efficiently in Gray code order. 
Indeed, for $A_2^n(1^{\ell-1} 0)$ with $\ell\geq 3$ it is enough to generate
(efficiently) the Gray code for $A_2^n(01^{\ell-1})$ (see Theorem \ref{the:general_max})
and then reverse each generated word; and for  $A_2^n(0^{\ell-1} 1)$ 
it is enough to generate the Gray code for $A_2^n(1^{\ell-1} 0)$ as previously, 
then complement each symbol in each word. The following scheme describes this method 
(see the example in Table \ref{ex_A2_4}):
$$A_2^n(01^{\ell-1})\xrightarrow{\text{   Reverse    
}}A_2^n(1^{\ell-1}0)\xrightarrow{\text{Complement}}A_2^n(0^{\ell-1}1).$$

Finally, notice that the generating order ($\prec$ or
$\triangleleft$ in our case) does not affect the efficiency of the generating algorithm,
which can obviously be modified to produce same set of factor avoiding words in lexicographical 
order. A {\tt C} implementation of our generating
algorithm is on the web site of the last author \cite{Vaj_web}.

\section{Conclusions}
We introduce two order relations on the set of length $n$ $q$-ary words,
and show that the set of words avoiding any from among the $q^\ell$ factors of length
$\ell\geq 2$, except $\ell-1$ or $\ell$ of them according to the parity of $q$, when listed in the appropriate
order is an (at most) $3$-Gray code.
For each of the excepted factors we give a simple transformation 
which allows to eventually obtain similar Gray codes. Finally, an efficient generating algorithm
for the derived Gray codes is given.

\clearpage
\section*{Appendix}
\begin{table}[!h]
\centering
\begin{tabular}{| c | c | c | c | c | c |}
\hline
0 0 0 0 		 & 0 1 \textbf{2} 2  & 1 0 \textbf{1 0}  & 1 \textbf{2} 2 0  &  2 2 \textbf{2} 2 & 2 0 \textbf{1} 2 \\
0 0 0 \textbf{1}  & 0 1 2 \textbf{1}  & 1 0 1 \textbf{1}  & 1 2 2 \textbf{1}  &  2 2 2 \textbf{1} & 2 0 1 \textbf{1} \\
0 0 0 \textbf{2}  & 0 1 2 \textbf{0}  & 1 0 1 \textbf{2}  & 1 2 2 \textbf{2}  &  2 2 2 \textbf{0} & 2 0 1 \textbf{0} \\
0 0 \textbf{1 0}  & 0 \textbf{2} 2 0  & 1 0 \textbf{2} 2  & 1 2 \textbf{1} 2  &  2 \textbf{1} 2 0 & 2 0 \textbf{0 2} \\
0 0 1 \textbf{1}  & 0 2 2 \textbf{1}  & 1 0 2 \textbf{1}  & 1 2 1 \textbf{1}  &  2 1 2 \textbf{1} & 2 0 0 \textbf{1} \\
0 0 1 \textbf{2}  & 0 2 2 \textbf{2}  & 1 0 2 \textbf{0}  & 1 2 1 \textbf{0}  &  2 1 2 \textbf{2} & 2 0 0 \textbf{0} \\
0 0 \textbf{2} 2  & 0 2 \textbf{1} 2  & 1 \textbf{1 0} 0  & 1 2 \textbf{0 2}  &  2 1 \textbf{1} 2 &  \\
0 0 2 \textbf{1}  & 0 2 1 \textbf{1}  & 1 1 0 \textbf{1}  & 1 2 0 \textbf{1}  &  2 1 1 \textbf{1} &  \\
0 0 2 \textbf{0}  & 0 2 1 \textbf{0}  & 1 1 0 \textbf{2}  & 1 2 0 \textbf{0}  &  2 1 1 \textbf{0} &  \\
0 \textbf{1 0} 0  & 0 2 \textbf{0 2}  & 1 1 \textbf{1 0}  & \textbf{2} 2 0 0  &  2 1 \textbf{0 2} &  \\
0 1 0 \textbf{1}  & 0 2 0 \textbf{1}  & 1 1 1 \textbf{1}  & 2 2 0 \textbf{1}  &  2 1 0 \textbf{1} &  \\
0 1 0 \textbf{2}  & 0 2 0 \textbf{0}  & 1 1 1 \textbf{2}  & 2 2 0 \textbf{2}  &  2 1 0 \textbf{0} &  \\
0 1 \textbf{1 0}  & \textbf{1 0} 0 0  & 1 1 \textbf{2} 2  & 2 2 \textbf{1 0}  &  2 \textbf{0 2} 0 &  \\
0 1 1 \textbf{1}  & 1 0 0 \textbf{1}  & 1 1 2 \textbf{1}  & 2 2 1 \textbf{1}  &  2 0 2 \textbf{1} &  \\
0 1 1 \textbf{2}  & 1 0 0 \textbf{2}  & 1 1 2 \textbf{0}  & 2 2 1 \textbf{2}  &  2 0 2 \textbf{2} &  \\

\hline
\end{tabular}
\caption{
\label{ex_A3_4}
The set $A_3^4$ listed in $\triangleleft$ 
order, inducing a 2-Gray code. The list is columnwise and the changed symbols are in bold.}
\end{table}

\begin{table}[!h]
\centering
\begin{tabular}{| c | c | c |}
\hline
$A_2^4(011)$ & $A_2^4(110)$ & $A_2^4(001)$\\
(a)&(b)&(c)\\
\hline\hline
0 0 0 0 & 0 0 0 0  & 1 1 1 1      \\
0 0 0 \textbf{1} & \textbf{1} 0 0 0  & \textbf{0} 1 1 1        \\
0 0 \textbf{1 0} & \textbf{0 1} 0 0   & \textbf{1 0} 1 1       \\
0 \textbf{1 0 1} & \textbf{1 0 1} 0  &   \textbf{0 1 0} 1   \\
0 1 0 \textbf{0} & \textbf{0} 0 1 0   &  \textbf{1} 1 0 1     \\
\textbf{1} 1 0 0 & 0 0 1 \textbf{1}   &  1 1 0 \textbf{0}     \\
1 1 0 \textbf{1} & \textbf{1} 0 1 1   &  \textbf{0} 1 0 0    \\
1 1 \textbf{1} 1 & 1 \textbf{1} 1 1  &  0 \textbf{0} 0 0    \\
1 1 1 \textbf{0} & \textbf{0} 1 1 1   &  \textbf{1} 0 0 0      \\
1 \textbf{0} 1 0 & 0 1 \textbf{0} 1   &  1 0 \textbf{1} 0     \\
1 0 \textbf{0 1} & \textbf{1 0} 0 1   & \textbf{0 1} 1 0       \\
1 0 0 \textbf{0} & \textbf{0} 0 0 1    &  \textbf{1} 1 1 0     \\
\hline

\end{tabular}
\caption{
\label{ex_A2_4}
(a) The set $A_2^4(011)$ listed in $\prec$ order, inducing 3-adjacent 
Gray code;  (b) the reverse of the list in (a), giving Gray code for 
$A_2^4(110)$; (c) the complement of the list in (b), giving Gray code for 
$A_2^4(001)$. The changed symbols are in bold}
\end{table}

\end{document}